\documentclass[a4paper,11pt]{article}
\usepackage[latin1]{inputenc}
\usepackage{amsmath}
\usepackage{amsfonts}
\usepackage{amssymb}
\usepackage[brazil,english]{babel}
\usepackage[T1]{fontenc}
\usepackage{ae}
\usepackage{amsmath,amsthm,amssymb}
\usepackage{graphicx}
\usepackage{color}
\usepackage{subfigure}
\usepackage{float}
\usepackage{setspace}
\usepackage{pdfpages}
\usepackage{mathrsfs}
\usepackage{multicol}

\advance\oddsidemargin by -0.5cm  \advance\textwidth by 1cm

\theoremstyle{definition}
\newtheorem{Def}{Definition}[section]
\newtheorem{Obs}{Remark}[section]
\newtheorem{Prop}{Proposition}[section]
\newtheorem{Lemma}{Lemma}[section]
\newtheorem{Theo}{Theorem}[section]
\newtheorem{Ex}{Example}[section]

\newcommand{\cqd}{{{\hfill $\rule{2.5mm}{2.5mm}$}}\vspace{0.5cm}}

\begin{document}

\begin{center}

\vspace{1cm}

 {\LARGE {\bf Discrete Conley Index for Zero-dimensional Basic Sets\\[3mm]
}}

\end{center}

\vspace{0.2cm}

\begin{center}
\begin{multicols}{2}
{\large {Mariana Gesualdi Villapouca}
\footnote{Supported by Capes 
} \ \ \ \ \
\ \ \ \ \ \ \ \   Ketty A. de Rezende} \footnote{Partially supported by CNPq under grant 302592/2010-5 and by FAPESP under grant 2012/18780-0.
}
\end{multicols}

\vspace{0.3cm}
\begin{multicols}{2}
\textit{Instituto de Matem\'{a}tica e Estat\'{i}stica, \\ Universidade Federal Fluminense, \\ 24.020-140 - Niter\'{o}i - RJ, Brasil}

\textit{Departamento de Matem\'{a}tica, \\ Universidade Estadual de Campinas, \\
13.083-859- Campinas - SP, Brazil}
\end{multicols}
\end{center}

\vspace{0.3cm}

\begin{abstract}
A theorem is established where the computation of the discrete Conley index for zero dimensional basic sets is given with respect to the dynamical information contained in the associated structure matrices. A classification of the reduced homology Conley index of a zero dimensional basic set in terms of its Jordan real form is presented.
\end{abstract}

{\small

\noindent {\bf MSC2010 subject classification:} Primary: 37B30; 37B10; 37C05. Secondary: 37D15; 37D05; 37B35.
}

\section*{Introduction}

\par The algebraic and topological techniques in dynamical systems have played a significant role as can be verified, for example, in Morse Theory and more recently in Conley's Index Theory. This theory was introduced by C.Conley in the late seventies and has been further developed in a variety of contexts. Initially defined for continuous dynamical systems on compact sets, see \cite{Conley}, \cite{Salamon} and \cite{Salamon_2} and  later for discrete dynamical systems.

\par The Conley index was originally defined for isolated invariant sets $S$ of continuous flows and defined to be the homotopy type of the pointed space obtained from a filtration pair associated to $S$. It was then proven that the Conley index is independent of the choice of the filtration pair and is a homotopy invariant. This is precisely the difficulty in the discrete case since there are examples where the homotopy type of the filtration pair of $S$ may depend on the choice of the pair.

\par Motivated by this difficulty several definitions were presented for the discrete case. See \cite{Franks_Richeson}, \cite{Mrozek}, \cite{Robbin_Salamon} and \cite{Szymczak}. The version we adopt herein is due to  Richeson \cite{Richeson} and called the \textbf{reduced homological Conley index}. The reason behind this choice is twofold: firstly it is computationally easier to work with and secondly, it is the form used in the definition of connection matrix pairs.

\par This index will be defined more precisely in Section~\ref{background}. It is defined for an isolated invariant set $S$ of a continuous map $f$ as
$$Con_{*}(S) = (CH_{*}(S), \chi_{*}(S)),$$
where  $CH_{*}(S) =$ $\displaystyle{ \bigcap_{n > 0}((f_{P})_{*})^{n}(H_{*}(N_{L}, [L]))}$ and $\chi_{*}(S): CH_{*}(S) \rightarrow CH_{*}(S)$ is an isomorphism induced by  $(f_{P})_{*}$. The map  $(f_{P})_{*}$ is the homology induced map of the pointed space map  $f_{p}: N_{L} \rightarrow N_{L}$  associated to the  filtration pair $P = (N, L)$.

\par The Conley index has the crucial property of continuation which guarantees the invariance of the index under small $C^{0}$ perturbations of the system. See \cite{Franks_Richeson} and \cite{Richeson}. The continuation property together with the fact that fitted diffeomorphisms are $C^{0}$-dense on the set of diffeomorphisms on closed manifolds, see \cite{Shub_Sullivan}, motivated our investigation of the discrete Conley index for zero dimensional basic sets.

\par In this article, our main result Theorem \ref{Theo1} establishes a relation between a homological invariant, namely the Conley index, and a dynamical invariant obtained from the structure matrix of a zero dimensional basic set. 
We then proceed to characterize the isomorphism of the Conley index of a zero dimensional basic set by using the real Jordan form of its structure matrix, Theorem \ref{Theo2}. Furthermore, by using Theorem \ref{Theo1}, we obtain an alternative formulation for the results on the homology Zeta function and on Morse inequalities of  diffeomorphisms with hyperbolic chain recurrents sets on closed manifolds found in \cite{Franks}.

\section{Background} \label{background}

\subsection{Diffeomorphism with hyperbolic chain recurrent set and zero dimensional basic sets}

\par In this section we will present some basic definitions and theorems on the dynamical behaviour of diffeomorphisms $f: M \rightarrow M$ where $M$ is a closed manifold with metric $d(\cdot,\cdot)$.

\begin{Def}
Let $f: M \rightarrow M$ be a diffeomorphism. A point $x \in M$ is called \textbf{chain recurrent} if given $\varepsilon > 0$ there exists $x_{1}=x, x_{2}, \ldots, x_{n-1}, x_{n}=x$ (with $n = n(\varepsilon)$) such that 
$$d( f(x_{i}), x_{i+1}) < \varepsilon,  \hspace*{0.2cm} \forall 1 \leq i < n.$$
Define \textbf{chain recurrent set}, $\mathcal{R}(f)$, as the set of all chain recurrent points of $f$. 
\end{Def}

\par By Smale's spectral decomposition, Theorem  \ref{decompspectral}, one has that if  $\mathcal{R}(f)$ is hyperbolic it admits a basic set decomposition. 

\begin{Def}
Let $f: M \rightarrow M$ be a diffeomorphism. A compact $f$-invariant set  $\Lambda$ is said to have a \textbf{hyperbolic structure} provided that the tangent bundle of $M$ restricted to $\Lambda$ can be written as a Whitney sum $T_{\Lambda}M = E^{u} \oplus E^{s}$ of subbundles $Df$-invariant, and there are constants $C>0$, $\lambda \in (0, 1)$, such that
$$\Vert Df^{n}(v) \Vert \leq C \lambda^{n} \Vert v \Vert, \forall \hspace*{0.1cm} v \in E^{s}, n \geq 0$$
and
$$\Vert Df^{n}(v) \Vert \geq C^{-1} \lambda^{-n} \Vert v \Vert, \forall \hspace*{0.1cm} v \in E^{u}, n \geq 0$$
\end{Def}

\begin{Theo}[\textbf{Smale's Spectral Decomposition}, \cite{Smale}] \label{decompspectral}
Let $f: M \rightarrow M$ be a diffeomorphism such that its chain recurrent set $\mathcal{R}(f)$ has a  hyperbolic structure. Then $\mathcal{R}(f)$ is a finite disjoint union of compact invariant sets $\Omega_{1}, \Omega_{2}, \ldots, \Omega_{n}$ and each $\Omega_{i}$ contains an orbit of the system which is dense in $\Omega_{i}$.
\end{Theo}

\par The  $\Omega_{i}$'s of Theorem \ref{decompspectral} are called \textbf{basic sets} of $f$, and the fiber dimension of the bundle $E_{\Omega_{i}}^{u}$ is called the \textbf{index} of $\Omega_{i}$.

\par Given a zero dimensional basic set $\Omega$ of a diffeomorphism $f: M \rightarrow M$ with hyperbolic chain recurrent set, one can associate to the dynamics of $f$ in $\Omega$, a matrix with $0, 1$ and $-1$ entries, called the structure matrix.

\par In order to define the structure matrix one needs the relation between a zero dimensional basic set and its subshift of finite type. We refer the reader to \cite{Bowen}.

\begin{Def}
The \textbf{subshift of finite type} determined by a finite set $\cal S$ and the relation $\rightarrow$ is a  homeomorphism $\sigma: \Sigma \rightarrow \Sigma$, where $\displaystyle{\Sigma \subset \displaystyle{\prod_{- \infty}^{\infty}} \cal S}$ is defined as $\Sigma = \lbrace s = (\ldots, s_{-1}, s_{0}, s_{1}, \ldots) \hspace*{0.1cm} \vert \hspace*{0.1cm} s_{i} \rightarrow s_{i+1} \hspace*{0.2cm} \mbox{for all} \hspace*{0.2cm} i \rbrace$ and $\sigma(s) = s^{\prime}$ where $s_{i}^{\prime} = s_{i-1}$, so $\sigma$ shifts to the right.
\end{Def}

\par The subshift of finite type  $\sigma$ determined by $\cal S$ and $\rightarrow$ is called a \textbf{vertex shift} associated to the matrix $A$ if one numbers the elements of $\cal S$ from 1 to $n$ ($n$ is the cardinality of $\cal S$) and the matrix $A_{n \times n}$ is given as follows:
$$A_{ij} = \left\lbrace  \begin{array}{ll}       
1, & \mbox{if} \hspace*{0.2cm} s_{i} \rightarrow s_{j} \\
0, & \mbox{otherwise}
\end{array} \right.$$

\par Bowen proves in \cite{Bowen} that if $\Omega$ is a zero dimensional basic set of a diffeomorphism $f$ with hyperbolic $\mathcal{R}(f)$, then $f$ restricted to $\Omega$ is topologically conjugate to a  vertex shift $\sigma(G): \Sigma_{G} \rightarrow \Sigma_{G}$ by  $h: \Sigma_{G} \rightarrow \Omega$ that has the property that the locally constant function
$$\Delta(x) = \left\lbrace   
\begin{array}{rll}
1 & \mbox{if} & Df_{x}: E^{u}_{x} \rightarrow E^{u}_{f(x)} \hspace{0.2cm} \mbox{preserves orientation}, \\
-1 & \mbox{if} & Df_{x}: E^{u}_{x} \rightarrow E^{u}_{f(x)} \hspace{0.2cm} \mbox{reverses orientation}
\end{array}
\right.$$
is constant in $h(C(k))$ for each $k$, where $C(k) = \lbrace c \in \Sigma_{G} \hspace*{0.1cm} \vert \hspace*{0.1cm} c_{0} = k  \rbrace$. Hence, one defines $\Delta_{k}$ as the value of $\Delta(x)$ in $h(C(k))$.

\begin{Def}
Define the \textbf{structure matrix} $A$ for a  basic set $\Omega$  by $A_{jk} = \Delta_{k} G_{jk}$.  
\end{Def}

\par Next we present a result of  Bowen and Franks \cite{Bowen_Franks} which associates a structure matrix to a homological invariant.

\begin{Def} \label{DefFiltration}
Let $f: M \rightarrow M$ be a diffeomorphism with hyperbolic chain recurrent set with basic sets  $\lbrace \Omega_{i} \rbrace_{i=0}^{n}$, then a \textbf{filtration associated to $f$} is a collection of submanifolds $M_{0} \subset M_{1} \subset \cdots \subset M_{n}=M$ such that
\begin{enumerate}
\item $f(M_{i}) \subset int \hspace{0.02cm}(M_{i})$;
\item $\displaystyle{\Omega_{i} = \bigcap_{n=-\infty}^{\infty} f^{n}(M_{i} \setminus M_{i-1})}$.
\end{enumerate}
\end{Def}

\par The existence of a filtration follows easily from the existence of a smooth Lyapunov function  $\phi: M \rightarrow \mathbb{R}$ for the diffeomorphism $f$, \cite{Franks}.

\begin{Def}
Let $h: V \rightarrow V$ and $h^{\prime}: V^{\prime} \rightarrow V^{\prime}$ be two endomorphisms where $V$ and $V^{\prime}$ are vector spaces. 
\begin{enumerate}
\item Define the \textbf{nonnilpotent part} $h^{+}$ of $h$ as the map induced by $h$ in the quocient space $\dfrac{V}{V_{0}}$ where $V_{0} = \lbrace  v \in V \hspace*{0.1cm} \vert \hspace*{0.1cm} h^{k}(v) = 0 \hspace{0.2cm} \mbox{for some} \hspace{0.2cm} k > 0 \rbrace$.
\item One says that $h$ and $h^{\prime}$ are \textbf{conjugate} if there exists an isomorphism $$\Psi: V \rightarrow V^{\prime} \hspace{0.1cm}\mbox{such that}\hspace{0.1cm}  h^{\prime}\circ \Psi = \Psi \circ h.$$
\end{enumerate}
\end{Def}

\begin{Theo}[\cite{Bowen_Franks}] \label{BowenFranks2}
Suppose that the diffeomorphism $f: M \rightarrow M$ has hyperbolic chain recurrent set, $\lbrace M_{j} \rbrace$ is a filtration associated to $f$ and, $\Omega_{i}$ is a zero dimensional basic set of index $u$. If $A$ is a structure matrix $n \times n$ for $\Omega_{i}$ and $\mathbb{F}$ is a field, then
\begin{enumerate}
\item The nonnilpotent part $A^{+}$ of $A: \mathbb{F}^{n} \rightarrow \mathbb{F}^{n}$ is conjugate to the nonnilpotent part $f_{*u}^{+}$ de $f_{*u}: H_{u}(M_{i}, M_{i-1}; \mathbb{F}) \rightarrow H_{u}(M_{i}, M_{i-1}; \mathbb{F})$, and
\item the map $f_{*k}: H_{k}(M_{i}, M_{i-1}; \mathbb{F}) \rightarrow H_{k}(M_{i}, M_{i-1}; \mathbb{F})$ is nilpotent for all $k \neq u$.
\end{enumerate}
\end{Theo}

\subsection{Reduced homological Conley index}

\par Let $U$ be an open subset of a locally compact metric space $X$ and $f: U \rightarrow X$ a continuous map.

\par Define the \textbf{maximal invariant subset} of $N$,  $Inv \hspace{0.02cm}(N)$ as the set of points $x \in N$ such that there exists a complete orbit  $\lbrace  x_{n}\rbrace_{- \infty}^{\infty} \subset N$ with $x_{0} = x$ and $f(x_{n}) =  x_{n+1}$ for all $n$. One says that a set $S$ is an \textbf{isolated invariant set} if there exists an isolating neighborhood $N$ (i.e., a compact $N$ with $Inv \hspace{0.02cm}(N) \subset int(N)$) such that $S = Inv \hspace{0.02cm}(N)$.

\par Given an isolated invariant set $S$ there exists a pair of compact spaces $(N, L)$ with $L \subset N$ contained in the interior of the domain of $f$ with the following properties:
\begin{enumerate}
\item $ \overline{N \setminus L}$ is an isolating neighborhood of $S$,
\item $L$ is a neighborhood of $N^{-} = \lbrace x \in N \hspace*{0.1cm} \vert \hspace*{0.1cm} f(x) \notin int \hspace{0.02cm} (N)  \rbrace$ in $N$, 
\item $f(L) \cap \overline{N \setminus L} = \emptyset$
\end{enumerate}
See \cite{Franks_Richeson}. This compact pair is called \textbf{filtration pair}.

\par Let $P=(N,L)$ be a  filtration pair of an isolated invariant set $S$ of a continuous map  $f$. The map $f$ induces a continuous  base-point preserving map $f_{P}: N_{L} \rightarrow N_{L}$ with the property $[L] \subset int \hspace{0.02cm} f^{-1}([L])$, where $N_{L}$ is the quotient space $N/L$ and $[L]$ is the collapsed set that is taken as base point(see \cite{Franks_Richeson}). This map is the \textbf{pointed space map associated to $P$}. 

\begin{Def} \label{RedIndexRicheson}
Define \textbf{reduced homological Conley index} of $S$ as
$$Con_{*}(S) = (CH_{*}(S), \chi_{*}(S)),$$
where $\displaystyle{CH_{*}(S) = \bigcap_{n > 0}((f_{P})_{*})^{n}(H_{*}(N_{L}, [L]))}$ and $\chi_{*}(S): CH_{*}(S) \rightarrow CH_{*}(S)$ is the automorphism induced by $(f_{P})_{*}: H_{*}(N_{L}, [L]) \rightarrow  H_{*}(N_{L}, [L]) $ which is the homology induced map of the pointed space map  $f_{P}: N_{L} \rightarrow N_{L}$.
\end{Def}

\par Henceforth, whenever we mention the Conley index we mean the reduced homological Conley index  defined above.

\section{The Conley Index and the Structure Matrix}
 
\par In this section we prove the relation between the Conley index isomorphism of a zero dimensional basic set with the nonnilpotent part of its structure matrix.

\begin{Lemma} \label{lemma1}
Let $M$ be a closed manifold. Suppose $f: M \rightarrow M$ is a diffeomorphism with hyperbolic chain recurrent set, $\lbrace \Omega_{i} \rbrace_{i=0}^{n}$  are basic sets of $f$ and $\lbrace M_{i} \rbrace_{i=0}^{n}$ is an associated filtration. Then

\begin{enumerate}
\item $P_{i} = (M_{i}, M_{i-1})$  is a filtration pair for $\Omega_{i}$, $ \forall i=0, \ldots, n$. 

\item  $f_{*k}: H_{k}(M_{i}, M_{i-1}) \rightarrow H_{k}(M_{i}, M_{i-1})$ and $(f_{P_{i}})_{*k}: H_{k}(M_{i}/M_{i-1}, [M_{i-1}] ) \rightarrow H_{k}(M_{i}/M_{i-1}, [M_{i-1}] )$ are conjugate for all $k$, where $f_{P_{i}}$ is the pointed space map associated to $P_{i} = (M_{i}, M_{i-1})$. 
\end{enumerate}
\end{Lemma}

\begin{demo} 1. Let $\phi: M \rightarrow \mathbb{R}$ be a smooth Lyapunov function such that $M_{i} = \phi^{-1}((-\infty, c_{i} ])$ for all $i=0, \ldots, n$. Then $M_{i}$ is compact and
$$\displaystyle{Inv \hspace{0.02cm} (\overline{M_{i} \setminus M_{i-1}}) = \bigcap_{n \in \mathbb{Z}} f^{n} (\overline{M_{i} \setminus M_{i-1}}) = \Omega_{i}},$$
$M_{i}^{-} =  \lbrace  x \in M_{i} \vert f(x) \notin int \hspace{0.02cm} M_{i}  \rbrace   = \emptyset \subset M_{i-1}$, and $f(M_{i-1}) \cap \overline{M_{i} \setminus M_{i-1}} = \emptyset$, i.e., $(M_{i}, M_{i-1})$ is a filtration pair for $\Omega_{i}$, $ \forall i= 0, \ldots, n$.

2. Since the diagram
$$\begin{array}{lcr}
N & \stackrel{f}{\longrightarrow} & N \\
\downarrow^{q} & & \downarrow^{q} \\
\dfrac{N}{L} & \stackrel{f_{P}}{\longrightarrow} & \dfrac{N}{L}
\end{array}$$
is comutative ($q: N \rightarrow \frac{N}{L}$ is a quotient map), it follows that $q \circ f = f_{P} \circ q$. Therefore, $q_{*} \circ f_{*} = (f_{P})_{*} \circ q_{*}$, i.e., $f_{*}$ is conjugate to $(f_{P})_{*}$.
\cqd
\end{demo}

\begin{Theo} \label{Theo1}
Let $M$ be a closed manifold and $f: M \rightarrow M$ be a diffeomorphism. Suppose $f$ has a hyperbolic chain recurrent set and $\Omega_{i}$ is a zero dimensional basic set of index $u$. If $A$ is an $n \times n$ structure matrix for $\Omega_{i}$, and $\mathbb{F}$ is a field, then
\begin{enumerate}
\item $\chi_{u}(\Omega_{i})$ is conjugate to the nonnilpotent part $A^{+}$ of $A: \mathbb{F}^{n} \rightarrow \mathbb{F}^{n}$;

\item $Con_{k}(\Omega_{i}) = (0, 0)$ for all $k \neq u$.
\end{enumerate}
where $Con_{*}(\Omega_{i}) = (CH_{*}(\Omega_{i}), \chi_{*}(\Omega_{i}))$ is the Conley index of $\Omega_{i}$. 
\end{Theo}

\begin{demo}Let $\lbrace M_{j} \rbrace$ be a filtration associated to $f$. It is well known that $A^{+}$ is conjugate to the nonnilpotent part  $f_{*u}^{+}$ of $f_{*u}: H_{u}(M_{i}, M_{i-1}) \rightarrow H_{u}(M_{i}, M_{i-1})$ (see Theorem \ref{BowenFranks2}). By Lemma \ref{lemma1} we have that
$$f_{*k}: H_{k}(M_{i}, M_{i-1}) \rightarrow H_{k}(M_{i}, M_{i-1})$$
and
$$(f_{P})_{*k}: H_{k}(M_{i}/M_{i-1}, [M_{i-1}] ) \rightarrow H_{k}(M_{i}/M_{i-1}, [M_{i-1}] )$$
are conjugate for all $k$ and consequently $f_{*u}^{+}$ and $(f_{P})_{*u}^{+}$ are also conjugate.

\par We will show that $(f_{P})_{*u}^{+} = \chi_{u}(\Omega_{i}) $. Since, $(H_{u}(M_{i}/M_{i-1}, [M_{i-1}]))_{0}$ is the set
$$\lbrace v \in H_{u}(M_{i}/M_{i-1}, [M_{i-1}])\hspace*{0.1cm} \vert \hspace*{0.1cm} (f_{P})_{*u}^{k}(v) = 0 \hspace{0.2cm} \mbox{for some} \hspace{0.2cm} k > 0   \rbrace$$
then $(H_{u}(M_{i}/M_{i-1}, [M_{i-1}]))_{0} = g Ker ((f_{P})_{*u})$. Hence,  
$$\dfrac{H_{u}(M_{i}/M_{i-1}, [M_{i-1}])}{(H_{u}(M_{i}/M_{i-1}, [M_{i-1}]))_{0}}     = \dfrac{H_{u}(M_{i}/M_{i-1}, [M_{i-1}])}{g Ker ((f_{P})_{*u})} = (gIm(f_{P})_{*u}) = CH^{*}(\Omega_{i}).$$ 

\par Therefore, $\chi_{u}(\Omega_{i}) = (f_{P})_{*u}^{+}$ and this proves item 1.

\par It is well known that the map $f_{*k}$ is nilpotent for all $k \neq u$ (see Theorem \ref{BowenFranks2}). Hence, $(f_{P})_{*k}$ is also nilpotent for all $k \neq u$. Thus,
$$CH_{k}(\Omega_{i}) =  \dfrac{H_{k}(M_{i}/M_{i-1}, [M_{i-1}])}{g Ker ((f_{P})_{*k})} \approx 0, \forall k \neq u,$$
which implies that $\chi_{k}(\Omega_{i}) = 0$ for all $k \neq u$.
\cqd
\end{demo}

\begin{Obs}
An immediate consequence of our Theorem \ref{Theo1} and the Lemma 6.2 of \cite{Franks} is that if $\mbox{dim} \hspace*{0.1cm} W^{u}(\Omega_{i}) < k$ or $\mbox{dim} \hspace*{0.1cm} W^{s}(\Omega_{i}) < (\mbox{dim} \hspace*{0.1cm} M) - k$, then $Con_{k}(\Omega_{i}) = (0, 0)$. 
\end{Obs}

\begin{Obs}
Note that any Smale diffeomorphism on a closed manifold satisfies Theorem \ref{Theo1}. Moreover, if this diffeomorphism is Morse-Smale we have that $\chi_{u}(\Omega_{i})$ is conjugate to 
$$A = \left(
\begin{array}{lcrlc}
0 & 1 & 0 & \cdots & 0 \\
0 & 0 & 1 & \cdots & 0\\
\vdots & \vdots & \vdots & \ddots & \vdots \\
0 & 0 & 0 & \cdots & 1 \\
\pm 1 & 0 & 0 & \cdots & 0 \\
\end{array}
\right)_{m \times m}$$
\end{Obs}

\begin{Ex} \label{ExHorseshoe}
\par We will compute the Conley index of a Smale Horseshoe, Figure \ref{Horseshoe}, embedded in $S^{2}$ using Theorem \ref{Theo1}.

\begin{figure}[H]
\centering
\includegraphics[scale=0.8]{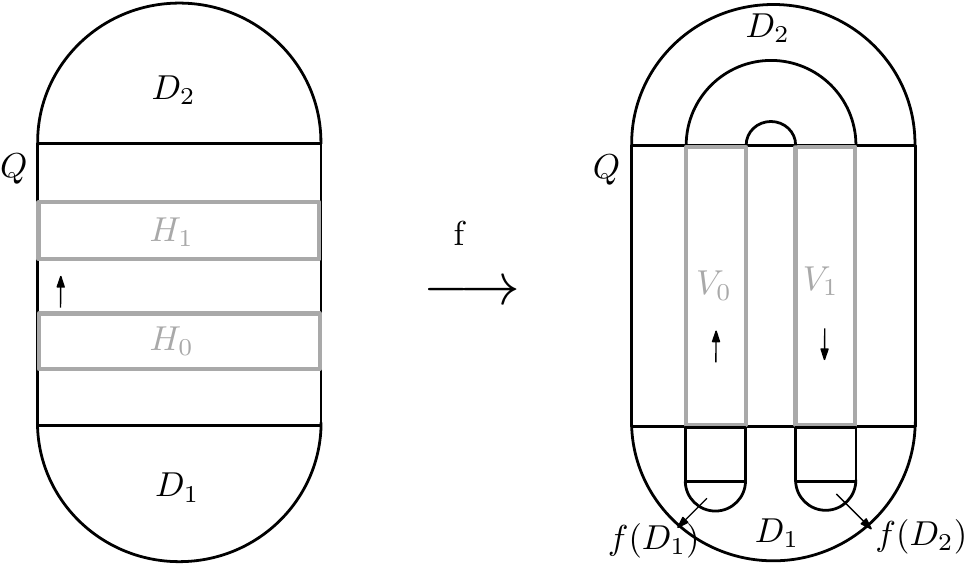}
\caption{Smale's Horseshoe} \label{Horseshoe}
\end{figure}

\par A structure matrix of the basic set $\displaystyle{\Lambda = \bigcap_{n \in \mathbb{Z}} f^{n}(Q)}$ of index 1 is given by
$$\left( \begin{array}{cc}  1 & -1 \\ 1 & -1 \end{array} \right)$$
with respect to the set of handles $H= H_{0} \cup H_{1}$.

Observe that $A^{+}=0$, since $A^{2} = 0$. Hence, the Conley index of Smale's Horseshoe is
$$Con_{q}(\Lambda) = (0,0), \forall q. $$

\end{Ex}

\par Observe that the advantage of this approach is that we do not make use of the filtration pair nor the pointed space map in order to compute the Conley index. We make use solely of the dynamical information contained in the structure matrix of Smale's Horseshoe.

\section{The Conley Index and the Real Jordan Form of a Structure Matrix}

\par In this section, we show that the isomorphism $\chi^{*}$, of the Conley index of a zero dimensional basic set of a diffeomorphism with hyperbolic chain recurrent set, is conjugate to the nonnilpotent part of an associated structure matrix. From a practical point of view, working with the structure matrix in its real Jordan form makes the  computation of the nilpotent and nonnilpotent part of the matrix effortless.
The information of the nonnilpotent part of this matrix is easily retrieved by a quick analysis of the non zero eigenvalues and the nilpotent part by the zero eigenvalues.

\par We prove the following theorem:

\begin{Theo} \label{Theo2}
Let $M$ be a closed manifold and $f: M \rightarrow M$ a diffeomorphism with a hyperbolic chain recurrent set and $\Omega$ a zero dimensional basic set of index $u$. Let $A$ be an $n \times n$ structure matrix for $\Omega$ and $\lambda_{1}, \ldots, \lambda_{r}$ their eigenvalues pairwise disjoint, then $\chi_{u}(\Omega)$ is conjugate to the real Jordan form of matrix $A$ after removing the Jordan blocks corresponding to zero eigenvalue. In other words, $\chi_{u}(\Omega)$ is conjugate to 
\begin{enumerate}
\item $J = \left( \begin{array}{ccccc} 
\lambda_{1} &              &             &           & 0        \\
            &  \lambda_{2} &             &           &           \\
            &              & \lambda_{3} &           &            \\
            &              &             &  \ddots   &             \\
 0          &              &             &           & \lambda_{n}
\end{array}  \right)$ if $r=n$ and $\lambda_{i} \neq  0, \forall i = 1, \ldots, n$.

\item $J = \left( \begin{array}{cccccc} 
\lambda_{1} &              &             &           &        & 0        \\
           &  \ddots &             &           &         &  \\
           &              & \lambda_{r_{0}-1} &           &    &        \\
                        &              &  &       \lambda_{r_{0}+1}    &        &    \\
            &              &             &     &  \ddots       &    \\
 0          &              &             &           &     & \lambda_{n}
\end{array}  \right)$ if $r=n$ and $\lambda_{r_{0}} = 0$.

\item $J = \left( \begin{array}{ccccc} 
J_{\lambda_{1}} &              &             &           & 0        \\
            & J_{\lambda_{2}} &             &           &           \\
            &              & J_{\lambda_{3}} &           &            \\
            &              &             &  \ddots   &             \\
 0          &              &             &           & J_{\lambda_{r}}
\end{array}  \right)$ if $r < n$ and $\lambda_{i} \neq  0, \forall i = 1, \ldots, r$.

\item $J = \left( \begin{array}{cccccc} 
J_{\lambda_{1}} &              &             &           &        & 0        \\
            &  \ddots &             &           &         &  \\
            &              & J_{\lambda_{r_{0}-1}} &           &    &        \\
                        &              &  &       J_{\lambda_{r_{0}+1}}    &        &    \\
            &              &             &     &  \ddots       &    \\
 0          &              &             &           &     & J_{\lambda_{r}}
\end{array}  \right)$  if $r< n$ and $\lambda_{r_{0}} = 0$.
\end{enumerate}
where for each  $i = 1, \ldots, r$ we have that $s_{i}$ is the geometric multiplicity of the eigenvalue $\lambda_{i}$, $n_{i}$ is the algebraic multiplicity of the eigenvalue $\lambda_{i}$ and $J_{\lambda_{i}}$ has the form
$$J_{\lambda_{i}} = \left(  \begin{array}{cccc}      
J_{\lambda_{i}, 1}  &                    &           &                          \\
                    & J_{\lambda_{i}, 2} &           &                           \\
                    &                    &  \ddots   &                            \\
                    &                    &           &  J_{\lambda_{i}, s_{i}}
\end{array}    \right)_{n_{i} \times n_{i}}$$
where each block $J_{\lambda_{i}, j}$ is a real Jordan block associated to the eigenvalue $\lambda_{i}$ which has one of the following forms:
$$\left(  \begin{array}{cccccc} 
\lambda_{i}       &      1     &       0        & \cdots     &     0     &   0    \\
0             & \lambda_{i}    &       1        &  \cdots    &     0     &   0     \\
\vdots        &   \vdots   &       \vdots   &   \ddots   &   \vdots  &  \vdots  \\
0             &      0     &       0        &   \cdots   &  \lambda_{i}  &   1      \\
0             &      0     &       0        &   \cdots   &  0        &   \lambda_{i}  
\end{array} \right)$$
if $\lambda_{i}$ is a real eigenvalue or
$$\left(  \begin{array}{cccccc} 
D      &     I     &       0        & \cdots     &     0     &   0    \\
0             & D   &       I        &  \cdots    &     0     &   0     \\
\vdots        &   \vdots   &       \vdots   &   \ddots   &   \vdots  &  \vdots  \\
0             &      0     &       0        &   \cdots   &  D  &   I      \\
0             &      0     &       0        &   \cdots   &  0        &   D 
\end{array} \right) $$
if $\lambda_{i} = a + b i$ is a complex eigenvalue, where $\displaystyle{ D = \left( \begin{array}{cc} a & -b \\ b & a \end{array} \right) \hspace{0.2cm} \mbox{and}   \hspace{0.2cm}  I = \left( \begin{array}{cc} 1 & 0 \\ 0 & 1  \end{array}  \right) }$.
\end{Theo}

\begin{demo} Let 
$$ J = \left( \begin{array}{cccccccccc} 
J_{\lambda_{1},1} &      &                    &                 &      &                        &      &                   &      &      0  \\       
                  & \ddots &                  &                 &      &                        &      &                   &       &       \\
                  &        & J_{\lambda_{1}, s_{1}}&             &      &                        &       &                  &       &       \\
                  &        &                    & J_{\lambda_{2}, 1}&   &                       &       &                  &       &       \\
                 &         &                                    &  &       \ddots               &       &                  &       &      \\
                 &         &                    &               &        &     J_{\lambda_{2}, s_{2}} &  &                 &       &      \\
                  &        &                    &               &        &                            & \ddots &           &       &      \\
                  &        &                    &               &        &                            &       & J_{\lambda_{r}, 1} &   &   \\  
                  &        &                    &               &        &                            &       &            &        \ddots \\
 0                 &        &                    &               &        &                            &       &            &       &       J_{\lambda_{r},s_{r}}
 \end{array} \right) $$
be the real Jordan form of $A$, where each $s_{i}$ is the geometric multiplicity of the eigenvalue $\lambda_{i}$ and each $J_{\lambda_{i}, j}$  is the real Jordan block associated to the eigenvalue $\lambda_{i}$.

\par Denote by $n_{i, j}$ the size of the real Jordan block $J_{\lambda_{i}, j}$ with $j = 1, \ldots, s_{i}$ and $i = 1, \ldots, r$. By the construction of the real Jordan form, we have that $n_{i,1} \geqslant n_{i, 2} \geqslant \cdots \geqslant n_{i, s_{i}}$. Since, $\displaystyle{n = \sum_{i=1}^{r} \sum_{j=1}^{s_{i}} n_{i,j}}$, then 
$$\mathbb{F}^{n} = [V_{\lambda_{1}, 1} \oplus \cdots \oplus V_{\lambda_{1}, s_{1}} ] \oplus [V_{\lambda_{2}, 1} \oplus \cdots \oplus V_{\lambda_{2}, s_{2}} ] \oplus \cdots \oplus [V_{\lambda_{r}, 1} \oplus \cdots \oplus V_{\lambda_{r}, s_{r}} ]$$
where
$$\begin{array}{cccc}
\displaystyle{ V_{\lambda_{1},1} = \bigoplus_{1}^{n_{1,1}} \mathbb{F},}  & \displaystyle{V_{\lambda_{1}, 2} = \bigoplus_{1}^{n_{1,2}} \mathbb{F},} & \cdots & \displaystyle{V_{\lambda_{1}, s_{1}} = \bigoplus_{1}^{n_{1,s_{1}}} \mathbb{F}} \\

\displaystyle{V_{\lambda_{2},1} = \bigoplus_{1}^{n_{2,1}} \mathbb{F},}  &\displaystyle{ V_{\lambda_{2}, 2} = \bigoplus_{1}^{n_{2,2}} \mathbb{F}}, & \cdots & \displaystyle{V_{\lambda_{2}, s_{2}} = \bigoplus_{1}^{n_{2,s_{2}}} \mathbb{F}} \\
  \vdots   &       \vdots    &     \vdots   &     \vdots   \\
     
\displaystyle{ V_{\lambda_{r},1} = \bigoplus_{1}^{n_{r,1}} \mathbb{F},}  & \displaystyle{V_{\lambda_{r}, 2} = \bigoplus_{1}^{n_{r,2}} \mathbb{F},} & \cdots & \displaystyle{V_{\lambda_{r}, s_{r}} = \bigoplus_{1}^{n_{r,s_{r}}} \mathbb{F}}

\end{array}$$

\par Thus, 
$$Ker J^{k} = [Ker J_{\lambda_{1}, 1}^{k} \oplus \cdots \oplus Ker J_{\lambda_{1}, s_{1}}^{k}] \oplus  \cdots \oplus [Ker J_{\lambda_{r}, 1}^{k} \oplus \cdots \oplus Ker J_{\lambda_{r}, s_{r}}^{k}], \forall k > 0, $$
$$gKer J = [gKer J_{\lambda_{1}, 1} \oplus \cdots \oplus gKer J_{\lambda_{1}, s_{1}}] \oplus  \cdots \oplus [gKer J_{\lambda_{r}, 1} \oplus \cdots \oplus gKer J_{\lambda_{r}, s_{r}}] $$
and
$$gIm J = [gIm J_{\lambda_{1}, 1} \oplus \cdots \oplus gIm J_{\lambda_{1}, s_{1}}] \oplus  \cdots \oplus [gIm J_{\lambda_{r}, 1} \oplus \cdots \oplus gIm J_{\lambda_{r}, s_{r}}] $$

\par Therefore,
$$J^{+} = \left( \begin{array}{ccccccc}    
J_{\lambda_{1},1}^{+}  &     &     &      &     &     &    0  \\

 &   \ddots  &     &      &     &     &      \\

 &     &   J_{\lambda_{1}, s_{1}}^{+}  &      &     &     &      \\

 &     &     &  \ddots   &     &     &      \\

 &     &     &     &    J_{\lambda_{r}, 1}^{+}  &    &      \\
 
  &     &     &     &      &  \ddots  &      \\
  
 0   &     &     &     &      &    &    J_{\lambda_{r}, s_{r}}^{+}  \\
\end{array} \right)$$

\par Hence, in order to calculate the nonnilpotent part, $J^{+}$, of the real Jordan form $J$ of $A$, it suffices to compute the  nonnilpotent part of the real Jordan blocks that make up $J$. In other words, it suffices to compute the nonnilpotent part of the following  real Jordan blocks associated to an eigenvalue $\lambda$:
\begin{itemize}
\item $J(1) =  \left(  \begin{array}{cccccc} 
\lambda       &      1     &       0        & \cdots     &     0     &   0    \\
0             & \lambda    &       1        &  \cdots    &     0     &   0     \\
\vdots        &   \vdots   &       \vdots   &   \ddots   &   \vdots  &  \vdots  \\
0             &      0     &       0        &   \cdots   &  \lambda  &   1      \\
0             &      0     &       0        &   \cdots   &  0        &   \lambda  
\end{array} \right)_{p \times p}$ if $\lambda \neq 0$

\item $J(2) =  \left(  \begin{array}{cccccc} 
0      &      1     &       0        & \cdots     &     0     &   0    \\
0             & 0    &       1        &  \cdots    &     0     &   0     \\
\vdots        &   \vdots   &       \vdots   &   \ddots   &   \vdots  &  \vdots  \\
0             &      0     &       0        &   \cdots   &  0  &   1      \\
0             &      0     &       0        &   \cdots   &  0        &  0  
\end{array} \right)_{p \times p}$ if $\lambda = 0$

\item $J(3) =\left(  \begin{array}{cccccc} 
D      &     I     &       0        & \cdots     &     0     &   0    \\
0             & D   &       I        &  \cdots    &     0     &   0     \\
\vdots        &   \vdots   &       \vdots   &   \ddots   &   \vdots  &  \vdots  \\
0             &      0     &       0        &   \cdots   &  D  &   I      \\
0             &      0     &       0        &   \cdots   &  0        &   D 
\end{array} \right) $ if $\lambda = a + b i$ is a complex eigenvalue ($b \neq 0$), where $\displaystyle{ D = \left( \begin{array}{cc} a & -b \\ b & a \end{array} \right) \hspace{0.2cm} \mbox{and}   \hspace{0.2cm}  I = \left( \begin{array}{cc} 1 & 0 \\ 0 & 1  \end{array}  \right) }$.
\end{itemize}

\par Firstly, we show that $J(1)^{+} = J(1)$. Since $\lambda \neq 0$, we have that
$$J(1)^{k} = \left(
\begin{array}{ccccc}
\lambda^{k}    &  C_{1} \lambda^{k-1}   &  \cdots  &  C_{p-2} \lambda^{k - (p-2)}   & C_{p-1} \lambda^{k-(p-1)} \\
0              & \lambda^{k}            &  \cdots  &  C_{p-3} \lambda^{k - (p-3)}   & C_{p-2} \lambda^{k - (p-2)} \\
\vdots         &  \vdots                &   \ddots &   \vdots                       &  \vdots                     \\
0              & 0                      &  \cdots  &   \lambda^{k}                  & C_{1} \lambda^{k-1}         \\
0              & 0                      & \cdots   &   0                            &  \lambda^{k}
\end{array}
\right)_{p \times p}$$
where $\displaystyle{C_{i} = \left(   \begin{array}{c} k \\ i \end{array}  \right) = \dfrac{k (k-1) \cdots (k - (i-1))}{i!}, \forall i = 1, \ldots, p-1}$. If $k \geqslant p$, then $C_{i} \neq 0$, for all $i = 1, \ldots, p-1$ and if $k < p$, then there exists $k \leqslant i_{0} < p$ such that $C_{i} \neq 0$ if $i = 1, \ldots, i_{0}$ and $C_{i} = 0$ for $i \geqslant i_{0}$. Hence, $gKer J(1) = \lbrace  (0, \ldots, 0) \rbrace$ and thus we conclude that $J(1)^{+} = J(1)$. 
 
\par Note that the matrix $J(2)$ is nilpotent, consequently $gKer J(2) = \mathbb{F}^{p}$, which implies that $J(2)^{+} = 0$.

\par Finally, we will prove that $J(3)^{+} = J(3)$. Indeed, since
$$J(3)^{k} = \left(
\begin{array}{ccccc}
D^{k}    &  C_{1} D^{k-1}   &  \cdots  &  C_{p^{\prime}-2} D^{k - (p^{\prime}-2)}   & C_{p^{\prime}-1} D^{k-(p^{\prime}-1)} \\
0              & D^{k}            &  \cdots  &  C_{p^{\prime}-3} D^{k - (p^{\prime}-3)}   & C_{p^{\prime}-2} D^{k - (p^{\prime}-2)} \\
\vdots         &  \vdots                &   \ddots &   \vdots                       &  \vdots                     \\
0              & 0                      &  \cdots  &   D^{k}                  & C_{1} D^{k-1}         \\
0              & 0                      & \cdots   &   0                            &  D^{k}
\end{array}
\right)$$
where $p^{\prime} = \dfrac{p}{2}$ and $C_{i} = \left(  \begin{array}{c} k \\ i \end{array} \right)  = \dfrac{k(k-1) \cdots (k - (i-1))}{i!}, \forall i = 1, \ldots, $ $p^{\prime} - 1$, we have that $det (J(3)^{k}) = det(D^{k}) \cdots  det(D^{k}) \neq 0$. This follows because $det D^{k} \neq 0$ which is easily proved by induction on  $k$. In this way, the columns of $J(3)^{k}$ form a linearly independent set. Therefore, $gKer J(3) = \lbrace (0, \ldots, 0)  \rbrace$, which implies that $J(3)^{+} = J(3)$. 
\cqd
\end{demo}

\par In what follows we present an example which confirms the practicability of the computation of the Conley index using Theorem \ref{Theo2}, specially in the cases of larger strucuture matrices.

\begin{Ex}
\par Consider four rectangles $H_{1}$, $H_{2}$, $H_{3}$ and $H_{4}$ in $\mathbb{R}^{2}$ and the diffeomorphism $f: \mathbb{R}^{2} \rightarrow \mathbb{R}^{2}$ shown in Figure \ref{Ex4squares}, where $f(H_{1}) = V_{1}$, $f(H_{2}) = V_{2}$, $f(H_{3}) = V_{3}$ and $f(H_{4}) = V_{4}$.
\begin{figure} 
\centering
\includegraphics[scale=0.5]{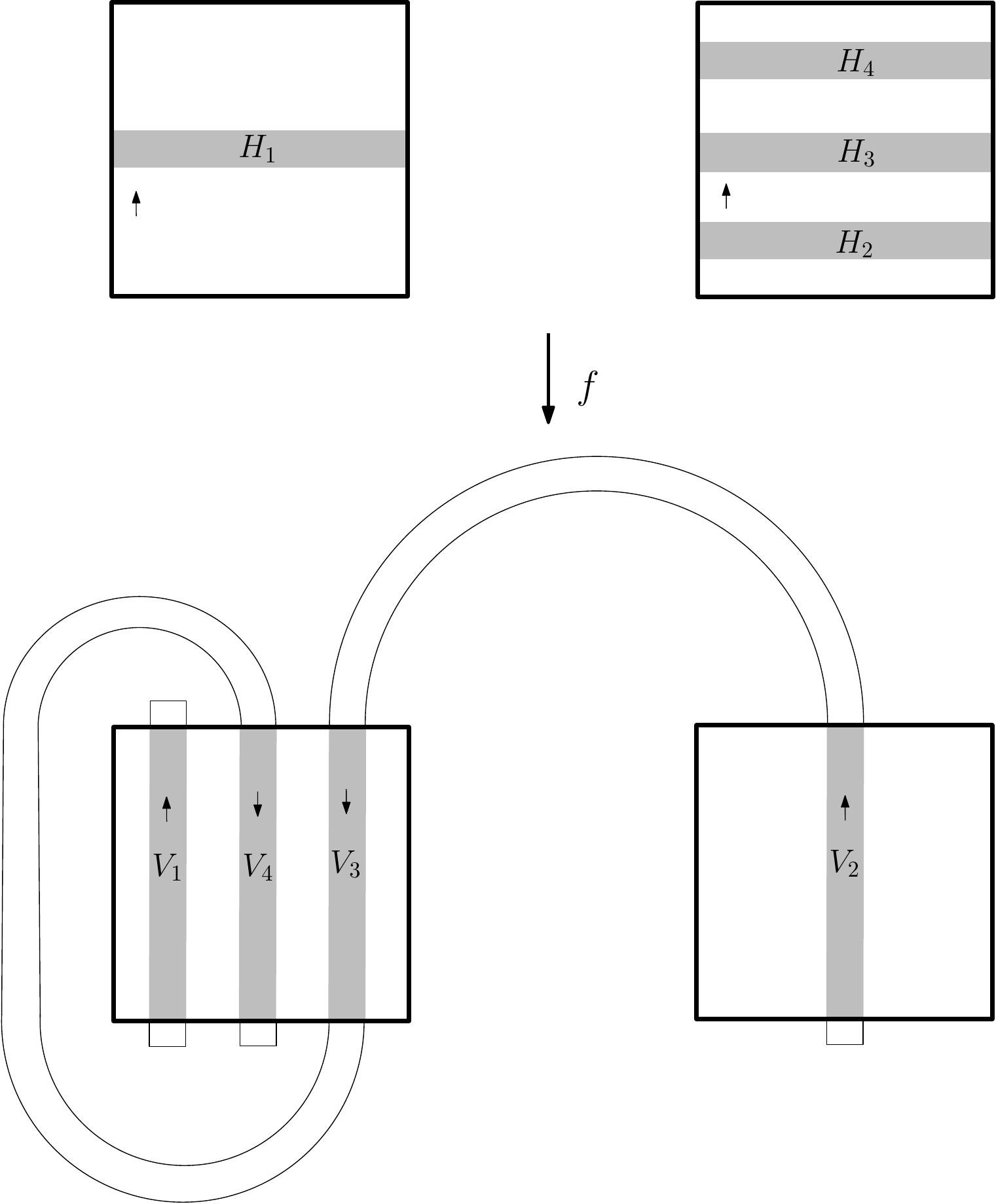}
\caption{fitted diffeomorphism}  \label{Ex4squares}
\end{figure}
\par A structure matrix of the basic set $\Lambda = \displaystyle{\bigcap_{n \in \mathbb{Z}} f^{n} (\bigcup_{i=1}^{4} H_{i})}$ of index 1 is
$$ A = \left(
\begin{array}{cccc}
1 & 0 & -1 & -1 \\
0 & 1 & 0 & 0 \\
0 & 1 & 0 & 0 \\
0 & 1 & 0 & 0
\end{array}
\right)$$
when we take as a set of handles $H = H_{1} \cup H_{2} \cup H_{3} \cup H_{4}$.
\par Since the real Jordan form of $A$ is
$$ \left(
\begin{array}{cccc}
0 & 0 & 0 & 0 \\
0 & 0 & 0 & 0 \\
0 & 0 & 1 & 1 \\
0 & 0 & 0 & 1
\end{array}
\right)$$
we have by Theorem \ref{Theo2}, that $\chi_{1}(\Lambda)$ is conjugate to $ \left(
\begin{array}{cc}
 1 & 1 \\ 
 0 & 1
\end{array}
\right)$ and therefore,
$$Con_{q}(\Lambda) = \left\lbrace \begin{array}{cc}
\left( \mathbb{R} \oplus \mathbb{R} , \left(   \begin{array}{cc} 1 & 1 \\ 0 & 1  \end{array} \right) \right), & q=1 \\ 
(0, 0), & q \neq 1
\end{array} \right.$$
\end{Ex}

\section{Conley Index, Homology Zeta Function and Morse Inequalities}

\par In this section we present the characterization of the homology Zeta function and of the Morse inequalities found in \cite{Franks}, using the isomorphism $\chi_{*}$ of the Conley index.

\begin{Prop} \label{Fczetaconley}
Let $M$ be a closed manifold. Suppose $f: M \rightarrow M$ is a diffeomorphism with hyperbolic chain recurrent set and $\lbrace \Omega_{i} \rbrace_{i=0}^{n}$ all the basic sets of $f$, then the homology Zeta function of $f\vert_{\Omega_{i}}$  is
$$Z_{i}(f) = \prod_{k = 0}^{dim M} det (I - \chi_{k}(\Omega_{i}) t)^{(-1)^{k+1}}$$
where $Con_{k}(\Omega_{i}) = (CH_{k}(\Omega_{i}), \chi_{k}(\Omega_{i}))$ is the Conley index with coefficients in $\mathbb{R}$ of the basic sets $\Omega_{i}$, $\forall i = 0, \ldots, l$.
\end{Prop}

\begin{demo} Let $\lbrace M_{i} \rbrace$ be a  filtration associated to $f$. We proved in Lemma \ref{lemma1}, that $$(f_{P})_{*k}: H_{k}(M_{i}/M_{i-1}, [M_{i-1}]; \mathbb{R}) \rightarrow H_{k}(M_{i}/M_{i-1}, [M_{i-1}]; \mathbb{R})$$ 
and 
$$f_{*k}: H_{k}(M_{i}, M_{i-1}; \mathbb{R}) \rightarrow H_{k}(M_{i}, M_{i-1}; \mathbb{R})$$ 
are conjugated by the isomorphism $q_{*k}: H_{k}(M_{i}, M_{i-1}; \mathbb{R}) \rightarrow  H_{k}(M_{i}/M_{i-1}, [M_{i-1}]; \mathbb{R}) $ induced by the quotient map $q: M_{i} \rightarrow \dfrac{M_{i}}{M_{i-1}}$. Thus,
$$\begin{array}{lll}
\mbox{det} \hspace{0.1cm}(I - (f_{P})_{*k} t) & = & \mbox{det} \hspace{0.1cm}(I - (q_{*k} \circ f_{*k} \circ q_{*k}^{-1}) t) \\

                        & = & \mbox{det} \hspace{0.1cm}(I - (q_{*k} \circ (f_{*k} t) \circ q_{*k}^{-1})) \\
                        
                        & = & \mbox{det} \hspace{0.1cm} (q_{*k} \circ (I - f_{*k} t) \circ q_{*k}^{-1}) \\
                        
                        & = & \mbox{det} \hspace{0.1cm}(q_{*k}) det(I - f_{*k} t) det(q_{*k})^{-1} \\
                        
                        & = & \mbox{det} \hspace{0.1cm}(I - f_{*k} t)
\end{array}$$

\par On the other hand, there exists a base $ H_{k}(M_{i}/M_{i-1}, [M_{i-1}]; \mathbb{R})$ such that
$$[(f_{P})_{*k}]_{\beta} =  \left( \begin{array}{cc}

(f_{P})_{*k}^{+}  &    0   \\

\ast &   (f_{P})_{*k}^{-}

\end{array} \right)$$
where $(f_{P})_{*k}^{-}$ is the nilpotent part of $(f_{P})_{*k}$. Hence, 
$$\mbox{det} \hspace{0.1cm}(I - (f_{P})_{*k}) = \mbox{det} \hspace{0.1cm}(I - (f_{P})_{*k}^{+}) \cdot \mbox{det} \hspace{0.1cm}(I - (f_{P})_{*k}^{-}) = \mbox{det} \hspace{0.1cm}(I - (f_{P})_{*k}^{+}),$$
follows since $(f_{P})_{*k}^{-} $ being nilpotent implies that $\mbox{det} \hspace{0.1cm}(I - (f_{P})_{*k}^{-}) = 1$.

\par Therefore,
$$\mbox{det} \hspace{0.1cm}(I - \chi_{k}(\Omega_{i}) t) = \mbox{det} \hspace{0.1cm}(I - (f_{P})_{*k}^{+} t) = \mbox{det} \hspace{0.1cm}(I - (f_{P})_{*k} t) = \mbox{det} \hspace{0.1cm}(I - f_{*k} t).$$ \cqd
\end{demo}

\par By using our Proposition \ref{Fczetaconley} and Theorem 6.5 of \cite{Franks}, we obtain the following proposition.

\begin{Prop} \label{IneqMC}
If $f: M \rightarrow M$ has hyperbolic chain recurrent set and their basic sets are homologically split at $q$ over $\mathbb{R}$, then there exists an integer polynomial $P(t)$ such that
$$P(t)^{(-1)^{q}} \prod_{u(i) \leqslant q} Z_{i}(f) = \prod_{k=0}^{q} det(I - \chi_{k}(M) t)^{(-1)^{k+1}}$$
where 
$$Z_{i}(f) = \prod_{k = 0}^{dim M} det (I - \chi_{k}(\Omega_{i}) t)^{(-1)^{k+1}}$$
and $Con_{k}(\Omega_{i}) = (CH_{k}(\Omega_{i}), \chi_{k}(\Omega_{i}))$  is the Conley index with coefficients in $\mathbb{R}$ of the basic set $\Omega_{i}$, $\forall i = 0, \ldots, l$.
\end{Prop}

\begin{demo} Since $P = (M, \emptyset)$ is a filtration pair of $M$, we have that $f_{*k}: H_{k}(M; \mathbb{R}) \rightarrow H_{k}(M; \mathbb{R})$ which is induced by $f$ and $(f_{P})_{*k}:H_{k}(M, \emptyset; \mathbb{R}) \rightarrow H_{k}(M, \emptyset; \mathbb{R})$ which is induced by  $f_{P}$ are isomorphic and therefore $\mbox{det} \hspace{0.1cm} (I - f_{*k} t) = \mbox{det} \hspace{0.1cm}(I - \chi_{k}(M) t)$. \cqd
\end{demo}

\par The construction of the following diffeomorphism of Example \ref{ExTorus} can be found in \cite{Franks}.

\begin{Ex} \label{ExTorus}
Consider the diffeomorphism $f: T^{2} \rightarrow T^{2}$ ilustrated in Figure \ref{SmaleTorus}.

\begin{figure}[h!] 
\begin{center}
\includegraphics[scale=1]{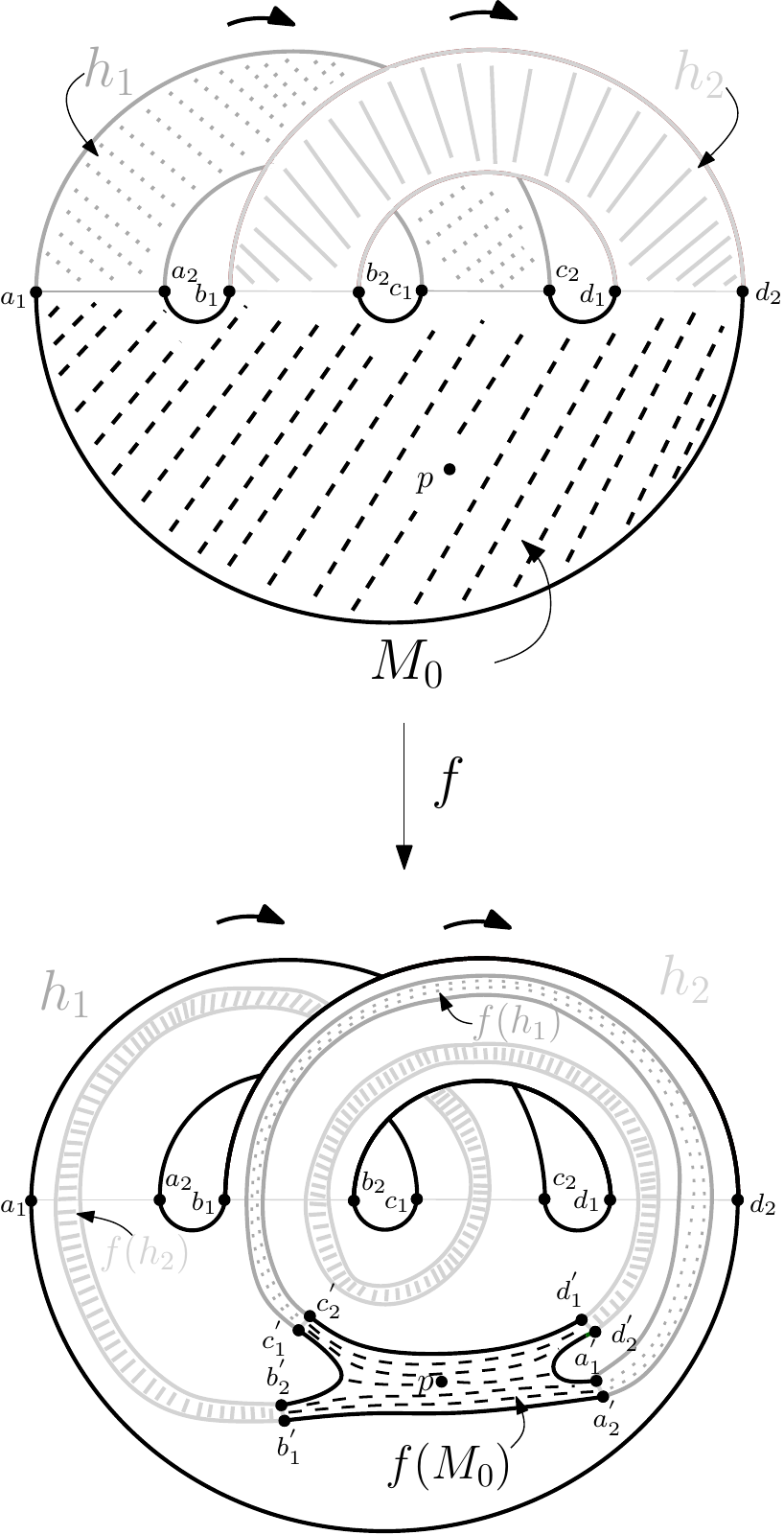} 
\caption{Smale diffeomorphism} \label{SmaleTorus}
\end{center}
\end{figure}

\par The diffeomorphism $f$ will have an attracting fixed point $p$, a repelling fixed point $\infty$ in the deleted disc $D_{ \infty }$ and is constructed in such a way that the two strips $H=h_{1} \cup h_{2}$ form a set of hyperbolic handles of $f$. Thus, $\mathcal{R}(f) = \lbrace p \rbrace \cup \Lambda \cup \lbrace \infty \rbrace$, where $\Lambda = \displaystyle{\bigcap_{n \in \mathbb{Z}} f^{n}(H)}$.

\par The structure matrices of $p$, $\Lambda$ and $\infty$ are respectively $A_{0} = (1)$, $A_{1} = \left(  \begin{array}{cc}  0 & 1 \\ -1 & 1    \end{array} \right)$ and $A_{2} = (1)$. Hence, by Theorem \ref{Theo1}, the Conley indices of the basic sets of $f$ are
$$\begin{array}{cc}
Con_{q}(p) = \left\lbrace \begin{array}{cc}
(\mathbb{R}, Id), & q=0 \\ 
(0, 0), & q \neq 0
\end{array} \right. ,  &

Con_{q}(\infty) = \left\lbrace \begin{array}{cc}
(\mathbb{R}, Id), & q=2 \\ 
(0, 0), & q \neq 2
\end{array} \right.

\end{array}$$
$$Con_{q}(\Lambda) = \left\lbrace \begin{array}{cc}
\left( \mathbb{R} \oplus \mathbb{R} , \left(   \begin{array}{cc} 0 & 1 \\ -1 & 1  \end{array} \right) \right), & q=1 \\ 
(0, 0), & q \neq 1
\end{array} \right.$$

\par Moreover, by Proposition \ref{Fczetaconley}, the homology Zeta functions associated to the basic sets are
$$ \begin{array}{lll}
Z_{p}(f) & = & \mbox{det} \hspace*{0.1cm} (I - \chi_{0}(p) t)^{-1} \hspace*{0.1cm} \mbox{det} \hspace*{0.1cm} (I - \chi_{1}(p) t) \hspace*{0.1cm} \mbox{det} \hspace*{0.1cm} (I - \chi_{2}(p) t)^{-1} \\
         & = & (1 - t)^{-1} \hspace*{0.1cm} (1) \hspace*{0.1cm} (1) = (1 - t)^{-1}
\end{array}$$

$$ \begin{array}{lll}
Z_{\Lambda}(f) & = & \mbox{det} \hspace*{0.1cm} (I - \chi_{0}(\Lambda) t)^{-1} \hspace*{0.1cm} \mbox{det} \hspace*{0.1cm} (I - \chi_{1}(\Lambda) t) \hspace*{0.1cm} \mbox{det} \hspace*{0.1cm} (I - \chi_{2}(\Lambda) t)^{-1} \\
         & = & (1) \hspace*{0.1cm} \mbox{det} \hspace*{0.1cm} \left( \begin{array}{cc} 1 & -t \\ t & (1-t)  \end{array} \right)        \hspace*{0.1cm} (1) = (1-t) + t^{2}
\end{array}$$

$$ \begin{array}{lll}
Z_{\infty}(f) & = & \mbox{det} \hspace*{0.1cm} (I - \chi_{0}(\infty) t)^{-1} \hspace*{0.1cm} \mbox{det} \hspace*{0.1cm} (I - \chi_{1}(\infty) t) \hspace*{0.1cm} \mbox{det} \hspace*{0.1cm} (I - \chi_{2}(\infty) t)^{-1} \\
         & = &  (1) \hspace*{0.1cm} (1) \hspace*{0.1cm} (1 - t)^{-1}  = (1 - t)^{-1}
\end{array}.$$

\par Now, since 
$$Con_{q}(T^{2}) = \left\lbrace \begin{array}{cc}
(\mathbb{R}, Id), & q=0 \\
\left( \mathbb{R} \oplus \mathbb{R} , \left(   \begin{array}{cc} 0 & 1 \\ -1 & 1  \end{array} \right) \right), & q=1 \\ 
(\mathbb{R}, Id), & q=2 \\
(0, 0), & \mbox{otherwise}
\end{array} \right.$$
and the basic sets of $f$ are homologically split at $1$ over $\mathbb{R}$ we have that, by Proposition \ref{IneqMC}, there exists an integer polynomial $P(t)$ such that
$$P(t)^{(-1)} \prod_{u(i) \leqslant 1} Z_{i}(f) = \prod_{k=0}^{1} det(I - \chi_{k}(M) t)^{(-1)^{k+1}}$$

\par Furthermore, we have that $P(t) = 1$, i.e.,
$$\prod_{u(i) \leqslant 1} Z_{i}(f) = \prod_{k=0}^{1} \mbox{det} \hspace{0.1cm}(I - \chi_{k}(M) t)^{(-1)^{k+1}},$$
since
$$\begin{array}{l}
\begin{array}{lll}
\displaystyle{\prod_{u(i) \leqslant 1} Z_{i}(f)} & = & Z_{0}(f) Z_{1}(f)  \\
                                                 & = & (1-t)^{-1} \hspace*{0.1cm} (t^{2} -t + 1)
                                                 
\end{array}    \\                                   

\begin{array}{lll}                                                 
\displaystyle{\prod_{k=0}^{1} \mbox{det} \hspace{0.1cm}(I - \chi_{k}(M) t)^{(-1)^{k+1}}} & = & \mbox{det} \hspace{0.1cm}(I - \chi_{0}(T^{2})t)^{-1} \mbox{det} \hspace{0.1cm}(I - \chi_{1}(T^{2})t) \\ 

    &  = &  (1-t)^{-1} \hspace*{0.1cm} (t^{2} -t + 1) 
\end{array}

\end{array} $$
\end{Ex}


\begin{thebibliography}{99}

\bibitem{Bowen} R. Bowen, \textit{Topological entropy and Axiom A,} Proc. Sympos. Pure Math., vol. 14, Amer. Math. Soc,. Providence, R.I., 1970,  23-41.


\bibitem{Bowen_Franks} R. Bowen and J. Franks, \textit{Homology for zero-dimensional nonwandering sets,} Annals of Mathematics \textbf{106} (1977), no. 1, 73-92. 


\bibitem{Conley} C. Conley, \textit{Isolated invariant sets and Morse index,} CBMS Regional Conference Series in Mathematics, n. 38, Amer. Math. Soc., Providence, R.I., 1978.


\bibitem{Franks} J. M. Franks, \textit{Homology and dynamical systems,} CBMS Regional Conference Series in Mathematics, n. 49, Amer. Math. Soc., Providence, R.I., 1982.


\bibitem{Franks_Richeson} J. M. Franks and D. S. Richeson, \textit{Shift equivalence and the Conley Index,} Transactions of the American Mathematical Society \textbf{352} (2000), no. 7, 3305-3322.


\bibitem{Mrozek} M. Mrozek, \textit{Leray functor and cohomological Conley index for discrete dynamical systems,} Transactions of the American Mathematical Society \textbf{318} (1990), no. 1, 149-178.


\bibitem{Richeson} D. S. Richeson, \textit{Connection matrix pairs,} Banach Center Publications \textbf{47} (1999), 219-232.


\bibitem{Robbin_Salamon} J. W. Robbin and D. Salamon,\textit{ Dynamical systems, shape theory and Conley index,} Ergodic Theory $\&$ Dynamical Systems \textbf{$8^{*}$} (1988), 375-393.


\bibitem{Salamon} D. Salamon, \textit{Connected simple systems and Conley index of isolated invariant sets,} Transactions of the American Mathematical Society \textbf{291} (1985), no. 1, 1-41.


\bibitem{Salamon_2} \underline{\makebox[1cm]{}} , \textit{Morse Theory, Conley index and Floer homology,} Bull. London Math. Soc. \textbf{22} (1990), 113-140.


\bibitem{Smale} S. Smale, \textit{Differentiable dynamical systems,} Bull. Amer. Math. Soc. \textbf{73} (1967), 797-817.


\bibitem{Shub_Sullivan} M. Shub and D. Sullivan, \textit{Homology and dynamical systems,} Topology \textbf{14} (1975), 109-132.


\bibitem{Szymczak} A. Szymczak, \textit{The Conley index for discrete semidynamical systems,} Topology and its Applications \textbf{66} (1995), 215-240. 

\end{thebibliography}
\end{document}